\documentclass{siamltex}
\usepackage[T1]{fontenc}
\usepackage[latin9]{inputenc}
\usepackage{float}
\usepackage{amssymb}
\usepackage{graphicx}

\makeatletter

\begin{document}

\title{Max-min and min-max approximation problems\\for normal matrices revisited}

\author{J\"org Liesen\footnotemark[1] and Petr Tich\'y\footnotemark[2]}

\footnotetext[1]{Institute of Mathematics, Technical University of Berlin,
Stra{\ss}e des 17. Juni 136, 10623 Berlin, Germany ({\tt liesen@math.tu-berlin.de}).}

\footnotetext[2]{Institute of Computer Science, Academy of Sciences of
the Czech Republic, Pod Vod\'arenskou v\v{e}\v{z}\'{\i} 2, 18207 Prague,
Czech Republic ({\tt tichy@cs.cas.cz}).
The work of this author was supported by the Grant Agency of the Czech Republic under grant
No. P201/13-06684~S, and by the project M100301201 of the institutional support
of the Academy of Sciences of the Czech Republic.}

\maketitle

\begin{center}
{\em In memory of Bernd Fischer}
\end{center}\medskip

\begin{abstract}
We give a new proof for an equality of certain max-min and min-max approximation 
problems involving normal matrices. The previously published proofs of this equality
apply tools from matrix theory, (analytic) optimization theory 
and constrained convex optimization. Our proof uses a classical characterization theorem 
from approximation theory and thus exploits the link between the two approximation problems 
with normal matrices on the one hand and approximation problems on compact sets in the 
complex plane on the other.
\end{abstract}



\section{Introduction}

Let $A$ be a real or complex square matrix, i.e., $A\in \mathbb{F}^{n\times n}$
with $\mathbb{F}=\mathbb{R}$ or $\mathbb{F}=\mathbb{C}$. Suppose that $f$ and
$\varphi_1,\dots,\varphi_k$ are given (scalar) functions so that
$f(A)\in \mathbb{F}^{n\times n}$ and $\varphi_1(A),\dots,\varphi_k(A)\in \mathbb{F}^{n\times n}$
are well defined matrix functions in the sense of~\cite[Definition~1.2]{B:Hi2008}.
(In the case $\mathbb{F}=\mathbb{R}$ this requires a subtle assumption which is
explicitly stated in (\ref{eqn:assumption}) below.)
Let ${\mathcal P}_k(\mathbb{F})$ denote the linear span of the functions $\varphi_1,\dots,\varphi_k$
with coefficients in ${\mathbb F}$, so that in particular $p(A)\in {\mathbb F}^{n\times n}$
for each polynomial $p=\alpha_1\varphi_1+\dots+\alpha_k\varphi_k \in {\mathcal P}_k(\mathbb{F})$.

With this notation, the optimality property of many useful methods of numerical
linear algebra can be formulated as an approximation problem of the form
\begin{equation} \label{eqn:minv}
\min_{p\in {\mathcal P}_k(\mathbb{F})}\|f(A)v-p(A)v\|,
\end{equation}
where $v\in {\mathbb F}^n$ is a given vector and $\|\cdot\|$ denotes the Euclidean norm
on $\mathbb{F}^n$. In~(\ref{eqn:minv}) we seek a best approximation (with respect to
the given norm) of the vector $f(A)v\in {\mathbb F}^n$ from the subspace of ${\mathbb F}^n$
spanned by the vectors $\varphi_1(A)v,\dots,\varphi_k(A)v$.
An example of such a method is the GMRES method~\cite{SaSc86}
for solving the linear algebraic problem $Ax=b$ with $A\in {\mathbb F}^{n\times n}$, $b\in {\mathbb F}^n$,
and the initial guess $x_0\in {\mathbb F}^n$. Its optimality
property is of the form (\ref{eqn:minv})
with $f(z)=1$, $\varphi_i(z)=z^i$ for $i=1,\dots,k$, and $v=b-Ax_0$.

If the given vector $v$ has unit norm, which usually can be assumed without loss of generality,
then an upper bound on (\ref{eqn:minv}) is given by
\begin{equation} \label{eqn:minmax}
\min_{p\in {\mathcal P}_k(\mathbb{F})}\|f(A)-p(A)\|,
\end{equation}
where $\|\cdot\|$ denotes the matrix norm associated with the Euclidean vector norm,
i.e., the matrix 2-norm or spectral norm on ${\mathbb F}^{n\times n}$.
In (\ref{eqn:minmax}) we seek a best approximation
(with respect to the given norm) of the matrix $f(A)\in {\mathbb F}^{n\times n}$
from the subspace of ${\mathbb F}^{n\times n}$ spanned by
the matrices $\varphi_1(A),\dots,\varphi_k(A)$. An example of this type is the
Chebyshev matrix approximation problem with $A\in {\mathbb F}^{n\times n}$,
$f(z)=z^k$, and $\varphi_i(z)=z^{i-1}$, $i=1,\dots,k$. This problem was introduced
in~\cite{GrTr1994} and later studied, for example, in~\cite{ToTr1998} and~\cite{FaLiTi2010}.

In order to analyse how close the upper bound (\ref{eqn:minmax}) can possibly be to the
quantity (\ref{eqn:minv}), one can maximize (\ref{eqn:minv}) over all unit norm vectors
$v\in {\mathbb F}^n$ and investigate the sharpness of the inequality
\begin{eqnarray}
\max_{v\in {\mathbb F}^n\atop \|v\|=1}\min_{p\in {\mathcal P}_k(\mathbb{F})}\|f(A)v-p(A)v\|
&\leq& \min_{p\in {\mathcal P}_k(\mathbb{F})}\|f(A)-p(A)\|  \label{eqn:ineq2}\\
&=&\min_{p\in {\mathcal P}_k(\mathbb{F})}\max_{v\in {\mathbb F}^n\atop \|v\|=1}\|f(A)v-p(A)v\|.\nonumber
\end{eqnarray}
From analyses of the GMRES method it is known that the inequality (\ref{eqn:ineq2}) can
be strict. For example, certain nonnormal matrices $A\in {\mathbb R}^{4\times 4}$ were 
constructed in~\cite{FaJoKnMa1996,To1997}, for which (\ref{eqn:ineq2}) is strict with $k=3$,
$f(z)=1$, and $\varphi_i(z)=z^i$, $i=1,2,3$.
More recently, nonnormal matrices $A\in{\mathbb R}^{2n\times 2n}$, $n\geq 2$, were derived
in~\cite{FaLiTi2013}, for which the inequality (\ref{eqn:ineq2}) is strict 
for all $k=3,\dots, 2n-1$, $f(z)=1$, and $\varphi_i(z)=z^i$, $i=1,\dots, k$.

On the other hand, the following result is well known.

\medskip
\begin{theorem}\label{thm:main}
Under the assumptions made in the first paragraph of the Introduction, if  
$A\in{\mathbb F}^{n\times n}$ is normal, then equality holds in (\ref{eqn:ineq2}).
\end{theorem}

\medskip
At least three different proofs of this theorem or variants of it can be found in the literature.
Greenbaum and Gurvits proved it for ${\mathbb F}={\mathbb R}$ using mostly
methods from matrix theory; see~\cite[Section~2]{GrGu1994} as well as Section~\ref{sec:final}
below for their formulation of the result. Using (analytic) methods of optimization
theory, Joubert proved the equality for the case of the GMRES method with $f(z)=1$, 
$\varphi_i(z)=z^i$, $i=1,\dots,k$, and he distinguished the cases ${\mathbb F}={\mathbb R}$ 
and ${\mathbb F}={\mathbb C}$; see~\cite[Theorem~4]{Jo1994}. Finally, Bellalij, Saad, and Sadok 
also considered the GMRES case with  ${\mathbb F}={\mathbb C}$, and they applied methods 
from constrained convex optimization; see~\cite[Theorem~2.1]{BeSaSa2008}.

In this paper we present yet another proof of Theorem~\ref{thm:main}, which is rather
simple because it fully exploits the link between matrix approximation problems for
normal matrices and scalar approximation problems in the complex plane. We observe
that when formulating the matrix approximation problems in (\ref{eqn:ineq2}) in
terms of scalar approximation problems, the proof of Theorem~\ref{thm:main} reduces
to a straightforward application of a well-known characterization theorem of best
approximation in the complex plane. While the proof of the theorem for ${\mathbb F}={\mathbb C}$
can be accomplished in just a few lines, the case ${\mathbb F}={\mathbb R}$ contains
some technical details that require additional attention. 

The characterization theorem from approximation theory we use in this paper
and some of its variants have been stated and applied also in other publications in
this context, and in particular in~\cite[Theorem~5.1]{BeSaSa2008}. To our knowledge
the theorem has, however, not been used to give a simple and direct proof of Theorem~\ref{thm:main}.

\medskip
{\em Personal note.} We have written this paper in memory of our colleague Bernd
Fischer, who passed away on July 15, 2013. Bernd's achievements in the analysis
of iterative methods for linear algebraic systems using results of approximation
theory, including his nowadays classical monograph~\cite{B:Fi1996}, continue to
inspire us in our own work. One of Bernd's last publications in this area
(before following other scientific interests), written jointly with
Franz Peherstorfer (1950--2009) and published 2001 in ETNA~\cite{FiPe2001},
also is based on a variant of the characterization theorem we apply in this paper.

\section{Characterization theorem and proof of Theorem~\ref{thm:main}}\label{sec:classical}
To formulate the characterization theorem of best approximation in the complex plane
we follow the treatment of Rivlin and Shapiro~\cite{RiSh1961} that has been 
summarized in Lorentz' book~\cite[Chapter~2]{B:Lo1986}.

Let $\Gamma$ be a compact subset of ${\mathbb F}$, and let $C(\Gamma)$ denote the set of continuous
functions on $\Gamma$. If $\Gamma$ consists of finitely many single points (which is the
case of interest in this paper), then $g\in C(\Gamma)$ means that the function $g$ has
a well defined (finite) value at each point of $\Gamma$. For $g\in C(\Gamma)$ we denote
the maximum norm on $\Gamma$ by
\[
\|g\|_{\Gamma}\equiv \max_{z\in\Gamma}|g(z)|.
\]

Now let $f\in C(\Gamma)$ and $\varphi_1,\dots,\varphi_k\in C(\Gamma)$ be given functions
with values in ${\mathbb F}$, where either ${\mathbb F}={\mathbb R}$
or ${\mathbb F}={\mathbb C}$. As above, let ${\mathcal P}_k(\mathbb{F})$ denote
the linear span of the functions $\varphi_1,\dots,\varphi_k$ with coefficients in ${\mathbb F}$.
For $p\in {\mathcal P}_k({\mathbb F})$, define
\[
\Gamma(p)\equiv\{z\in\Gamma:\ |f(z)-p(z)|=\|f-p\|_\Gamma\}.
\]

A function $p_*=\alpha_1\varphi_1+\dots+\alpha_k\varphi_k\in {\mathcal P}_k(\mathbb{F})$ is called
a {\em polynomial of best approximation} for $f$ on $\Gamma$ when
\begin{equation}\label{eqn:pba}
\|f-p_*\|_{\Gamma}\;=\;\min_{p\in {\mathcal P}_k(\mathbb{F})}\|f-p\|_\Gamma.
\end{equation}
Under the given assumptions such a polynomial of best approximation exists;
see, e.g.,~\cite[Theorem~1, p.~17]{B:Lo1986}.
The following well known result (see, e.g.,~\cite[Theorem~3, p.~22]{B:Lo1986}
or~\cite[pp. 672-674]{RiSh1961}) characterizes the polynomials of best approximation.

\medskip
\begin{theorem}\label{thm:RiSh}
In the notation established above, the following two statements are equivalent:
\begin{enumerate}
\item The function $p_*\in\mathcal{P}_k(\mathbb{F})$ is a polynomial of best approximation for $f$ on $\Gamma$.
\item For the function $p_*\in\mathcal{P}_k(\mathbb{F})$ there exist
$\ell$ pairwise distinct points $\mu_{i}\in\Gamma(p_*)$, $i=1,\dots,\ell$,
where $1\leq\ell\leq k+1$ for ${\mathbb F}={\mathbb R}$ and $1\leq\ell\leq2k+1$ for ${\mathbb F}={\mathbb C}$,
and $\ell$ real numbers $\omega_{1},\dots,\omega_{\ell}>0$ with $\omega_{1}+\dots+\omega_{\ell}=1$,
such that
\begin{equation}\label{eq:minmaxgen}
\sum_{j=1}^{\ell}\omega_{j}\,[f(\mu_{j})-p_*(\mu_{j})]\overline{p(\mu_{j})}=0,\quad
\mbox{for all $p\in\mathcal{P}_k(\mathbb{F})$.}
\end{equation}
\end{enumerate}
\end{theorem}

\medskip
A well known geometric interpretation of the condition (\ref{eq:minmaxgen}) is that the origin
is contained in the convex hull of the points
$$\left\{\left([f(\mu)-p_*(\mu)]\overline{\varphi_1(\mu)},\dots,[f(\mu)-p_*(\mu)]\overline{\varphi_k(\mu)}\right)
\in {\mathbb F}^k\;:\;\mu\in\Gamma(p_*)\,\right\};$$
see, e.g.,~\cite[Equation~(5), p.~21]{B:Lo1986}. Here we will not use this interpretation, but
rewrite (\ref{eq:minmaxgen}) in
terms of an algebraic orthogonality condition involving vectors and matrices. Using that condition we
will be able to prove Theorem~\ref{thm:main} in a straightforward way. We will distinguish the cases of
complex and real normal matrices, because the real case contains some subtleties.

\subsection{Proof of Theorem~\ref{thm:main} for ${\mathbb F}={\mathbb C}$}
Let $A\in {\mathbb C}^{n\times n}$ be normal. Then $A$ is unitarily diagonalizable,
$A=Q\Lambda Q^H$ with $\Lambda={\rm diag}(\lambda_1,\dots,\lambda_n)$ and $QQ^H=Q^HQ=I_n$.
In the notation established above, let $\Gamma=\{\lambda_1,\dots,\lambda_n\}$ and
suppose that $p_*\in {\mathcal P}_k(\mathbb{C})$ is a polynomial of best approximation
for $f$ on $\Gamma$, so that statement 2. from Theorem~\ref{thm:RiSh} applies to $p_*$.
With this setting, the matrix approximation problem (\ref{eqn:minmax}) can be seen 
as the scalar best approximation problem (\ref{eqn:pba}), i.e.,
$$\min_{p\in {\mathcal P}_k(\mathbb{C})}\|f(A)-p(A)\| =
\min_{p\in {\mathcal P}_k(\mathbb{C})}\|f(\Lambda)-p(\Lambda)\| 
=\min_{p\in {\mathcal P}_k(\mathbb{C})} \|f-p\|_{\Gamma}\,.$$
Without loss of generality we may assume that the eigenvalues of $A$ are ordered so
that $\lambda_j=\mu_j$ for $j=1,\dots,\ell$. We denote
$$\delta\equiv \|f-p_*\|_\Gamma=|f(\lambda_j)-p_*(\lambda_j)|\quad
\mbox{for $j=1,\dots,\ell$.}$$
Next, we define the vector
\begin{equation}\label{eqn:vstar}
v_*\equiv Q\xi,\;\;\mbox{where}\;\;
\xi\equiv \left[\xi_1,\dots,\xi_\ell,0,\dots,0\right]^T\in {\mathbb C}^n,\;\;
|\xi_{j}|^{2}=\omega_{j}, \;j=1,\dots,\ell.
\end{equation}
Since $Q$ is unitary and $\omega_{1}+\dots+\omega_{\ell}=1$,
we have $\|v_*\|=1$.

The condition (\ref{eq:minmaxgen}) can be written as
\begin{eqnarray*}
0 & = & \sum_{j=1}^{\ell}|\xi_{j}|^{2}\overline{p(\lambda_{j})}\left[f(\lambda_{j})-p_*(\lambda_{j})\right]
= \xi^{H}p(\Lambda)^{H}\left[f(\Lambda)-p_*(\Lambda)\right]\xi\\
& = & v_*^{H}p(A)^{H}\left[f(A)-p_*(A)\right]v_*\,,
\quad\mbox{for all $p\in {\mathcal P}_k(\mathbb{C})$,}
\end{eqnarray*}
or, equivalently,
\begin{equation}\label{eqn:cond}
f(A)v_*-p_*(A)v_*\perp p(A)v_*\,,\quad \mbox{for all $p\in {\mathcal P}_k(\mathbb{C})$.}
\end{equation}
It is well known that this algebraic orthogonality condition with respect to the Euclidean
inner product is equivalent to the optimality condition
\begin{equation}\label{eqn:cond1}
\|f(A)v_*-p_*(A)v_*\|=\min_{p\in {\mathcal{P}}_k(\mathbb{C})}\|f(A)v_*-p(A)v_*\|;
\end{equation}
see, e.g.,~\cite[Theorem~2.3.2]{LieStrBook12}.

Using the previous relations we now obtain
\begin{eqnarray*}
\min_{p\in {\mathcal P}_k(\mathbb{C})}\|f(A)-p(A)\| &=&\delta=
\left(\sum_{j=1}^\ell |\xi_{j}|^{2}\delta^2\right)^{1/2}\\
&=&\left(\sum_{j=1}^\ell |\xi_{j}|^{2}\left|f(\lambda_{j})-p_*(\lambda_{j})\right|^2\right)^{1/2} \\
&=&\|\left[f(\Lambda)-p_*(\Lambda)\right]\xi\|\\
&=& \|Q\left[f(\Lambda)-p_*(\Lambda)\right]Q^H Q \xi \|\\
&=&\|f(A)v_*-p_*(A)v_*\| \\
&=&\min_{p\in {\mathcal{P}}_k(\mathbb{C})}\|f(A)v_*-p(A)v_*\|\\
&\leq& \max_{v\in {\mathbb C}^n\atop \|v\|=1}\min_{p\in {\mathcal P}_k(\mathbb{C})}\|f(A)v-p(A)v\|.
\end{eqnarray*}
This is just the reverse of the inequality (\ref{eqn:ineq2}) for ${\mathbb F}={\mathbb C}$,
and hence the proof of Theorem~\ref{thm:main} for ${\mathbb F}={\mathbb C}$ is complete.

\subsection{Proof of Theorem~\ref{thm:main} for ${\mathbb F}={\mathbb R}$}\label{sec:real_case}
Let $A\in {\mathbb R}^{n\times n}$ be normal. Since $A$ is real, its eigenvalues are either
real, or non-real but then occur in complex conjugate pairs. Therefore, the diagonal 
matrix $\Lambda={\rm diag}(\lambda_1,\dots,\lambda_n)$ and the unitary matrix $Q$ in the spectral decomposition~$A=Q\Lambda Q^H$ are in general complex.
Hence, the proof given in the previous section requires some modifications because the 
vector $v_*=Q\xi$ constructed in (\ref{eqn:vstar}) is then also complex in general, 
while for a real matrix $A$ the maximization in (\ref{eqn:ineq2}) is done over $v\in {\mathbb R}^n$.
Note that the proof presented in the previous section applies to a real and symmetric $A$, 
since then the matrices $Q$ and $\Lambda$ are real, giving a real vector $v_*=Q\xi$ in
(\ref{eqn:vstar}).

As above, let $\Gamma=\{\lambda_1,\dots,\lambda_n\}$. Since $A$ is real and normal, 
the set $\Gamma$ may contain non-real points, and for a general proof we must allow complex-valued
functions $f\in C(\Gamma)$ and $\varphi_1,\dots,\varphi_k\in C(\Gamma)$. This means that we must
work with Theorem~\ref{thm:RiSh} for ${\mathbb F}={\mathbb C}$, although $A$ is real. 
However, we will assume that for each eigenvalue $\lambda_j$ of $A$ the given functions $f$
and $\varphi_1,\dots,\varphi_k$ satisfy
\begin{equation}\label{eqn:assumption}
\overline{f(\lambda_j)}=f(\overline{\lambda}_j)\quad\mbox{and}\quad
\overline{\varphi_i(\lambda_j)}=\varphi_i(\overline{\lambda}_j),\;\; i=1,\dots,k.
\end{equation}
This is a natural assumption for real matrices $A$, since it guarantees that the matrices
$f(A)$ and $\varphi_1(A),\dots,\varphi_k(A)$ are real as well; 
see~\cite[Remark~1.9]{B:Hi2008}
(for analytic functions it is actually a necessary and sufficient condition; 
see~\cite[Theorem~1.18]{B:Hi2008}).

Now let $q_{*}=\sum_{i=1}^k\alpha_i\varphi_i\in\mathcal{P}_{k}(\mathbb{C})$ 
be a polynomial of best approximation for $f$ on~$\Gamma$. Then, for any eigenvalue 
$\lambda_j$ of $A$,
$$\big| f(\lambda_j) - \sum_{i=1}^k \alpha_i \varphi_i(\lambda_j) \big|=
\big| \overline{f(\lambda_j)} - \sum_{i=1}^k \overline{\alpha}_i \overline{\varphi_i(\lambda_j)} \big|=
\big| f(\overline{\lambda}_j) - \sum_{i=1}^k \overline{\alpha}_i \varphi_i(\overline{\lambda}_j) \big|.$$
Since both $\lambda_j$ and $\overline{\lambda}_j$ are elements of $\Gamma$, we see that also
$\overline{q}_{*}\equiv \sum_{i=1}^k\overline{\alpha}_i\varphi_i$ is a polynomial of 
best approximation for $f$ on~$\Gamma$. Denote
$$\delta\equiv \left\Vert f-q_{*}\right\Vert_\Gamma =\left\Vert f-\overline{q}_{*}\right\Vert_\Gamma,$$
then for any $0\leq\alpha\leq 1$ we obtain
\begin{eqnarray*}
\delta & \leq & \left\Vert f-\alpha q_{*}-(1-\alpha)\overline{q}_{*}\right\Vert _{\Gamma}=\left\Vert \alpha(f-q_{*})+(1-\alpha)(f-\overline{q}_{*})\right\Vert _{\Gamma}\\
 & \leq & \alpha\left\Vert f-q_{*}\right\Vert _{\Gamma}+
(1-\alpha)\left\Vert f-\overline{q}_{*}\right\Vert _{\Gamma}=\delta,
\end{eqnarray*}
which shows that any polynomial of the form $\alpha q_{*}+(1-\alpha)\overline{q}_{*}$,
$0\leq\alpha\leq1$, is also a polynomial of best approximation for $f$ on $\Gamma$. 
In particular, for $\alpha=\frac12$ we obtain the {\em real} polynomial
of best approximation 
$$p_{*}\equiv\frac{1}{2}(q_{*}+\overline{q}_{*}) \;\in\;\mathcal{P}_{k}(\mathbb{R}).$$
Using $p_*\in \mathcal{P}_k(\mathbb{R})$ and  (\ref{eqn:assumption}) we get 
$$|f(z)-p_*(z)| =
    |\overline{f(z)-p_*(z)}| =
    |f(\overline{z})-p_*(\overline{z})|,\quad\mbox{for all $z\in\Gamma$.}$$
Therefore, the set $\Gamma(p_*)$ of all points $z$ which satisfy
$|f(z)-p_*(z)| = \| f-p_*\|_\Gamma$ is symmetric with respect to the real axis, i.e.,
$z \in \Gamma(p_*)$ if and only if $\overline{z} \in \Gamma(p_*)$.

For simplicity of notation we denote
\[
\zeta_p(z) \equiv [f(z)-p_*(z)]\overline{p(z)}.
\]
In the definition of $\zeta_p(z)$ we indicate only its dependence on $p$ and $z$,
since $f$ is a given function and $p_*$ is fixed. If $p\in \mathcal{P}_k(\mathbb{R})$,
then the corresponding function $\zeta_p(z)$ satisfies 
$\zeta_p(z) = \zeta_p(\overline{z})$ for all $z\in\Gamma$.

\medskip
Now, Theorem~\ref{thm:RiSh} (with ${\mathbb F}={\mathbb C}$) implies the existence of a set
$$G_*\equiv \{\mu_1,\dots,\mu_\ell\}\subseteq \Gamma(p_*)\subseteq\Gamma,$$
and the existence of positive real numbers $\omega_1,\dots,\omega_\ell$ with $\sum_{j=1}^\ell \omega_j=1$,
such that
\begin{equation}\label{eqn:rel1}
\sum_{j=1}^{\ell}\omega_{j}\,\zeta_p(\mu_j)=0,\quad
\mbox{for all $p\in\mathcal{P}_k(\mathbb{R})$,}
\end{equation}
where we have used that $\mathcal{P}_k(\mathbb{R})\subset \mathcal{P}_k(\mathbb{C})$.
To define a convenient real vector $v$ similar to the construction leading to
(\ref{eqn:vstar}), we will ``symmetrize'' the condition (\ref{eqn:rel1}) with respect
to the real axis.

Taking complex conjugates in (\ref{eqn:rel1}), and using that
 $\zeta_p(z) = \zeta_p(\overline{z})$ for any $z\in\Gamma$, we obtain another relation of the form
\begin{equation}\label{eqn:rel2}
\sum_{j=1}^{\ell}\omega_{j}\,\zeta_p(\overline{\mu}_j)=0,\quad
\mbox{for all $p\in\mathcal{P}_k(\mathbb{R})$,}
\end{equation}
and, therefore,
\begin{equation}\label{eqn:rel3}
\frac12\sum_{j=1}^{\ell}\omega_{j}\,\zeta_p(\mu_j) + \frac12\sum_{j=1}^{\ell}\omega_{j}\,\zeta_p(\overline{\mu}_j)=0,\quad
\mbox{for all $p\in\mathcal{P}_k(\mathbb{R})$.}
\end{equation}
Here (\ref{eqn:rel3}) is the desired ``symmetrized'' condition. We now define
the set
\[
G_*^{\mathrm{sym}} \,\equiv\, \{\theta_1,\dots,\theta_m\}\,\equiv\, G_* \cup \overline{G}_*,
\]
Each $\theta_i\in G_*^{\mathrm{sym}}$ corresponds to some $\mu_j$ or $\overline{\mu}_j$,
and clearly $\ell\leq m\leq 2\ell$. (The exact value of $m$ is unimportant for our construction.) Writing
the condition (\ref{eqn:rel3}) as a single sum over all points from $G_*^{\mathrm{sym}}$, we get
\begin{equation}\label{eqn:rel4}
\sum_{i=1}^{m}\widetilde{\omega}_{i}\,\zeta_p(\theta_i)=0,\quad
\mbox{for all $p\in\mathcal{P}_k(\mathbb{R})$.}
\end{equation}
where the coefficients $\widetilde{\omega}_{i}$ are defined as follows:

If $\mu_j\in {\mathbb R}$, then $\zeta_p(\mu_j)$ appears in both sums in (\ref{eqn:rel3})
with the same coefficient $\omega_j/2$. Therefore, the term $\zeta_p(\theta_i)$ appears
in (\ref{eqn:rel4}) with the coefficient $\widetilde{\omega}_{i}=\omega_{j}$.

If $\mu_j \notin \mathbb{R}$ and $\overline{\mu}_j \notin G_*$, then $\zeta_p(\mu_j)$ appears
only in the left sum in (\ref{eqn:rel3}) with the coefficient $\omega_j/2$. Therefore, the term $\zeta_p(\mu_j)$
corresponds to a single term $\zeta_p(\theta_i)$ in (\ref{eqn:rel4}) with the
coefficient $\widetilde{\omega}_{i}=\omega_{j}/2$. Similarly, $\zeta_p(\overline{\mu}_j)$ appears only
in the right sum in (\ref{eqn:rel3}) with the coefficient $\omega_j/2$, and it corresponds to a single
term, say $\zeta_p(\theta_s)$, in (\ref{eqn:rel4}) with the coefficient $\widetilde{\omega}_{s}=\omega_{j}/2$.

If $\mu_j \notin \mathbb{R}$ and $\overline{\mu}_j \in G_*$, then $\overline{\mu}_j = \mu_s$ for
some index $s \neq j$, $1\leq s \leq \ell$. Therefore, the term $\zeta_p(\mu_j)$ appears in both sums in (\ref{eqn:rel3}),
in the left sum with the coefficient $\omega_j/2$, and in the right sum with the coefficient $\omega_s/2$. Hence,
$\zeta_p(\mu_j)$ corresponds to a single term $\zeta_p(\theta_i)$ in (\ref{eqn:rel4}) with the coefficient $\widetilde{\omega}_{i}=\omega_{j}/2+\omega_{s}/2$. Similarly, $\zeta_p(\overline{\mu}_j)$ corresponds to the
term $\zeta_p(\overline{\theta}_i)$ in (\ref{eqn:rel4}) with the coefficient equal to $\omega_{j}/2+\omega_{s}/2$.

One can easily check that $\widetilde{\omega}_{i}>0$ for $i=1,\dots,m$, and
\[\sum_{i=1}^m \widetilde{\omega}_{i} = 1.\]
Moreover, if
$\theta_j = \overline{\theta}_i$ for $j\neq i$, then
$\widetilde{\omega}_{j}=\widetilde{\omega}_{i}$.

Based on the relation (\ref{eqn:rel4}) we set
\begin{equation}\label{eqn:vstarreal}
v_*\equiv Q\xi,\;\;
\xi\equiv \left[\xi_1,\dots,\xi_n\right]^T\in {\mathbb R}^n,
\end{equation}
where  the $\xi_j$, $j=1,\dots,n$, are defined as follows:
If $\lambda_j \in G_*^{\mathrm{sym}}$, then there exits an index $i$ such that
$\lambda_j=\theta_i$ and we define $\xi_j \equiv \sqrt{\widetilde{\omega}_i}$.
If $\lambda_j \notin G_*^{\mathrm{sym}}$, we set $\xi_j = 0$.

It remains to  justify that the resulting vector $v_*$ is real.
If $\lambda_j\in {\mathbb R}$, then the corresponding eigenvector $q_j$ (i.e., the $j$th column of the
matrix $Q$) is real, and $\xi_j q_j$ is real. If $\lambda_j\notin {\mathbb R}$ and $\lambda_j \in G_*^{\mathrm{sym}}$,
then also $\overline{\lambda}_j\in G_*^{\mathrm{sym}}$, and $\overline{\lambda}_j=\lambda_i$ for
some $i\neq j$. The corresponding eigenvector is $q_i=\overline{q}_j$, and since $\xi_i=\xi_j$,
the linear combination $\xi_jq_j+\xi_iq_i=\xi_j(q_j+\overline{q}_j)$ is a real vector. Therefore,
the resulting vector $v_* = Q\xi$ is real.

Using (\ref{eqn:rel4}), analogously to the previous section, we get
$$ 0 = v_*^{T}p(A)^{T}\left[f(A)-p_*(A)\right]v_*\,,
\quad\mbox{for all $p\in {\mathcal P}_k(\mathbb{R})$,}$$
or, equivalently,
$$\|f(A)v_*-p_*(A)v_*\|=\min_{p\in {\mathcal{P}}_k(\mathbb{R})}\|f(A)v_*-p(A)v_*\|,$$
so that
\begin{eqnarray*}
\min_{p\in {\mathcal P}_k(\mathbb{R})}\|f(A)-p(A)\|
&=& \delta=\|f(A)v_*-p_*(A)v_*\| \\
&=&\min_{p\in {\mathcal{P}}_k(\mathbb{R})}\|f(A)v_*-p(A)v_*\|\\
&\leq& \max_{v\in {\mathbb R}^n\atop \|v\|=1}\min_{p\in {\mathcal P}_k(\mathbb{R})}\|f(A)v-p(A)v\|.
\end{eqnarray*}
This is just the reverse of the inequality (\ref{eqn:ineq2}) for ${\mathbb F}={\mathbb R}$,
and hence the proof of Theorem~\ref{thm:main} for ${\mathbb F}={\mathbb R}$ is complete.

\section{A different formulation}\label{sec:final}

Theorem~\ref{thm:main} can be easily rewritten as a statement about pairwise commuting normal matrices.
In the following we only discuss the complex case. The real case requires an analogous treatment as in
Section~\ref{sec:real_case}.

Let $A_0,A_1,\dots,A_k\in {\mathbb C}^{n\times n}$ be pairwise commuting normal matrices. Then
these matrices can be simultaneously unitarily diagonalized, i.e., there exists a unitary
matrix $U\in {\mathbb C}^{n\times n}$ so that
$$U^HA_iU=\Lambda_i={\rm diag}(\lambda_1^{(i)},\dots,\lambda_n^{(i)}),\quad i=0,1,\dots,k;$$
see, e.g.,~\cite[Theorem 2.5.5]{B:HoJo1990}.
Let $\Gamma\equiv \{\lambda_1,\dots,\lambda_n\}$
be an arbitrary set containing $n$ pairwise distinct complex numbers, and
let $A\equiv{\rm diag}(\lambda_1,\dots,\lambda_n)\in{\mathbb C}^{n\times n}$.
We now {\em define} the functions $f\in C(\Gamma)$ and $\varphi_1,\dots,\varphi_k\in C(\Gamma)$ 
to be any functions satisfying
$$f(\lambda_j)\equiv \lambda_j^{(0)},\quad \varphi_i(\lambda_j)\equiv \lambda_j^{(i)},\quad
j=1,\dots,n,\;\;i=1,\dots,k.$$
Then $f(A)=\Lambda_0$ and $\varphi_i(A)=\Lambda_i$ for $i=1,\dots,k$, 
so that Theorem~\ref{thm:main} implies
\begin{eqnarray*}
\max_{v\in\mathbb{C}^n\atop \|v\|=1} \min_{\alpha_1,\dots,\alpha_k\in {\mathbb C}}
\|A_0v-\sum_{i=1}^k \alpha_iA_iv\| &=&
\max_{v\in\mathbb{C}^n\atop \|v\|=1} \min_{\alpha_1,\dots,\alpha_k\in {\mathbb C}}
\|\Lambda_0v-\sum_{i=1}^k \alpha_i\Lambda_iv\| \\
&=&\max_{v\in\mathbb{C}^n\atop \|v\|=1} \min_{\alpha_1,\dots,\alpha_k\in {\mathbb C}}
\|f(A)v-\sum_{i=1}^k \alpha_i\varphi_i(A)v\| \\
&=&\min_{\alpha_1,\dots,\alpha_k\in {\mathbb C}}\|f(A)-\sum_{i=1}^k \alpha_i\varphi_i(A)\| \\
&=&\min_{\alpha_1,\dots,\alpha_k\in {\mathbb C}}\|A_0-\sum_{i=1}^k \alpha_iA_i\|.
\end{eqnarray*}
This equality is in fact the version of Theorem~\ref{thm:main} proven by Greenbaum and Gurvits
in~\cite[Theorem~2.3]{GrGu1994} for the case ${\mathbb F}={\mathbb R}$.

\def\polhk#1{\setbox0=\hbox{#1}{\ooalign{\hidewidth
  \lower1.5ex\hbox{`}\hidewidth\crcr\unhbox0}}}

%\bibliographystyle{siam}
%\bibliography{paper_normal}

\end{document}